\documentclass[amsfonts, 12pt]{article}
\usepackage{amssymb}

\newtheorem{theorem}{Theorem}[section]
\newtheorem{corollary}[theorem]{Corollary}
\newtheorem{example}[theorem]{Example}
\newtheorem{proposition}[theorem]{Proposition}
\newtheorem{lemma}[theorem]{Lemma}
\newtheorem{definition}[theorem]{Definition}
\newtheorem{remark}[theorem]{Remark}

\def\cB{\mathcal{B}}

\def\cD{\mathcal{D}}
\def\cE{\mathcal{E}}

\def\cF{\mathcal{F}}
\def\cL{\mathcal{L}}
\def\cM{\mathcal{M}}
\def\cN{\mathcal{N}}

\def\cP{\mathcal{P}}
\def\cS{\mathcal{S}}

\def\bC{\mathbb{C}}

\def\bR{\mathbb{R}}

\def\ve{\varepsilon}

\begin{document}

\title{Integration with respect to L\'evy colored noise,
with applications to SPDEs}

\author{Raluca M. Balan\footnote{University of Ottawa, Department of Mathematics and Statistics,
585 King Edward Avenue, Ottawa, ON, K1N 6N5, Canada. E-mail
address: rbalan@uottawa.ca} \footnote{Research supported by a
grant from the Natural Sciences and Engineering Research Council
of Canada.}}

\date{July 25, 2013}
\maketitle

\begin{abstract}
\noindent In this article, we introduce a L\'evy analogue of the spatially homogeneous Gaussian noise of \cite{dalang99}, and we construct a stochastic integral with respect to this noise. The spatial covariance of the noise is given by a tempered measure $\mu$ on $\bR^d$, whose density is given by $|h|^2$ for a complex-valued function $h$. Without assuming that the Fourier transform of $\mu$ is a non-negative function, we identify a large class of integrands with respect to this noise. As an application, we examine the linear stochastic heat and wave equations driven by this type of noise.
\end{abstract}



\section{Introduction}

There are two approaches for the study of stochastic partial differential equations (SPDEs) in the literature, known as the Da Prato-Zabczyk approach and the Walsh approach, initiated in the landmark references \cite{DZ92} and \cite{walsh86}, respectively. Depending on which approach one is using, the solution of an SPDE on $\bR_{+} \times \bR^d$ can be viewed as a process $(X_t)_{t \geq 0}$ with values in a suitable space of functions on $\bR^d$ (in the Da Prato-Zabczyk approach), or as a real-valued random field $\{u(t,x); t \geq 0,x \in \bR^d\}$ (in the Walsh approach). These two approaches have evolved independently and each has been fruitful in its own way. A direct comparison of the results obtained using the two approaches is not possible (see \cite{dalang-quer11} for some recent results in this direction, for equations with Gaussian noise).

In the recent years, a lot of attention has been given to the study of SPDEs with L\'evy noise, without Gaussian component. A comprehensive treatment can be found in the monograph \cite{PZ07}, using the Da Prato-Zabczyk approach. The goal of the present article is to introduce the basic tools which are necessary for initiating a similar study using the Walsh approach.

A basic model for a L\'evy noise without a Gaussian component can be defined using the same steps as in It\^o's classical construction of a L\'evy process. If the construction is done on $\bR_{+} \times \bR^d$, one obtains a space-time L\'evy white noise. (The details of this construction are given in Section \ref{white-noise-section} below.) This process is related to the ``impulsive cylindrical noise'' of \cite{PZ07}. Similarly to the Gaussian space-time white noise defined in \cite{walsh86}, the space of (deterministic) integrands with respect to the L\'evy white noise is $L^2(\bR_{+} \times \bR^d)$. Therefore, with this noise, even basic SPDEs (like the heat equation or the wave equation) will have random field solutions only in dimension $d=1$.

To avoid this problem, we introduce a {\em L\'evy colored noise} which can be viewed
as a counterpart of the spatially homogeneous Gaussian noise considered by Dalang in \cite{dalang99}. 
This noise is given by a process $X=\{X_t(\varphi);t \geq 0, \varphi \in \cS(\bR^d)\}$ defined via a representation based on the ``Fourier transform'' in space of the L\'evy white noise (see Definition \ref{def-colored-noise}), and has covariance
\begin{equation}
\label{cov-X}E[X_t(\varphi)X_s(\psi)]=(t \wedge s)\int_{\bR^d} \cF \varphi(\xi) \overline{\cF \psi(\xi)}\mu(d\xi),
\end{equation}
where $\cF \varphi$ is the Fourier transform of $\varphi$, and $\mu$ is a tempered measure on $\bR^d$ with density given by $|h|^2$, for a complex-valued function $h$.
Clearly, $X$ is spatially homogeneous, in the sense that for any $\varphi,\psi \in \cS(\bR^d)$ and $h \in \bR^d$,
$$E[X_t(\tau_h \varphi)X_t(\tau_h \psi)]=E[X_t(\varphi)X_t(\psi)],$$
where $\tau_h \varphi$ is the translation of $\varphi$ by $h$, i.e. $(\tau_h \varphi)(x)=\varphi(x+h)$ for all $x \in \bR^d$.

Under some additional assumptions (which are not needed in the present work), the L\'evy colored noise can also be constructed as an integral with respect to a compensated Poisson random measure on $\bR_{+} \times \bR^d \times (\bR \verb2\2 \{0\})$
(see Remark \ref{impulsive-remark} below). This construction has lead the authors of \cite{PZ06} to call it an ``impulsive colored noise''. A study of SPDEs with spatially homogeneous L\'evy noise (in particular, impulsive colored noise) can be found in Chapter 14 of \cite{PZ07}, using the Da Prato-Zabczyk approach.

In the present article, we develop a theory of stochastic integration with respect to $X$, using the same tools from Fourier analysis as in \cite{dalang99}. But unlike \cite{dalang99}, for this theory we do not require that the Fourier transform of $\mu$ be a non-negative function. Our main result (Theorem \ref{main}) identifies a large class of integrands with respect to $X$, which includes processes with values in the space $\cS'(\bR^d)$ of tempered distributions, and shows that for integrands in this class, the stochastic integral with respect to $X$
admits the same spectral representation as the process $X$ itself. (A similar result has been recently found in \cite{BGP12}, for the Gaussian noise.) This result allows us to study some linear SPDEs (like the heat or wave equations) with L\'evy colored noise, in any space dimension $d$.


This article is organized as follows. In Section \ref{white-noise-section} we  construct the L\'evy white noise and we examine its properties. In
Section \ref{Levy-section}, we give the definition of the L\'evy colored noise, we construct a stochastic integral with respect to this noise, and we identify a large class of integrands. In Section \ref{spde-section}, we give an application to the study of some SPDEs, like the stochastic heat equation and the stochastic wave equation.

We conclude the introduction with few words about the notation. We denote by
$\cB_b(\bR^d)$ the class of all bounded Borel sets in $\bR^d$ and by $|B|$ the Lebesgue measure of a set $B$ in $\bR^d$. We let
$L^2(\bR^d)$ be the set of all square-integrable functions on $\bR^d$, $\cD(\bR^d)$ be the set of all infinitely differentiable functions on $\bR^d$ with compact support, and $\cS(\bR^d)$ be the set of all infinitely differentiable functions on $\bR^d$ with rapid decrease. The analogue sets for complex-valued functions are denoted by $L_{\bC}^2(\bR^d)$, $\cD_{\bC}(\bR^d)$, respectively $\cS_{\bC}(\bR^d)$. Similar notations are used for the spaces $\bR^{d+1}$ and $\bR_{+} \times \bR^d$. We denote by
$\cS'(\bR^d)$ the class of tempered distributions on $\bR^d$. We let $\cF \varphi$ be the Fourier transform of a function $\varphi$ in $\cS(\bR^d)$ (or $L^2(\bR^d)$).

\section{The L\'evy white noise}
\label{white-noise-section}

In this section, we introduce the space-time L\'evy  white noise. This process plays an important role in the present article and can be viewed as an analogue of the space-time white noise introduced by Walsh in \cite{walsh86}.

We begin by generalizing to higher dimensions It\^o's construction of a classical L\'evy process. We refer the reader to Section 5.5 of \cite{resnick07} for an excellent pedagogical account of this construction.

Let $N=\sum_{i \geq 1}\delta_{(T_i,X_i,Z_i)}$ be a Poisson random measure on $\bR_{+} \times \bR^d \times (\bR \verb2\2 \{0\})$ defined on a probability space $(\Omega,\cF,P)$, with intensity measure $dtdx\nu(dz)$ where $\nu$ is a L\'evy measure on $\bR$, i.e. $\nu(\{0\})=0$ and
$$\int_{\bR} (1 \wedge |z|^2) \nu(dz)<\infty.$$ Let $(\varepsilon_j)_{j \geq 0}$ be a sequence of positive real numbers such that $\ve_j \to 0$ as $j \to \infty$ and $1=\ve_0>\ve_1>\ve_2>\ldots$. Let
$$\Gamma_j=\{z \in \bR; \ve_{j}<|z| \leq \ve_{j-1}\}, \ j\geq 1 \quad \mbox{and} \quad \Gamma_0=\{z \in \bR; |z|>1\}.$$

For any set $B \in \cB_b(\bR_{+} \times \bR^d)$, we define
$$L_j(B)=\int_{B \times \Gamma_j}z N(dt,dx,dz)=\sum_{(T_i,X_i) \in B}Z_i 1_{\{Z_i \in \Gamma_j\}}, \quad j \geq 0.$$

\begin{remark}
{\rm The variable $L_0(B)$ is finite since the sum above contains finitely many terms. To see this, we note that $E[N(B \times \Gamma_0)]=|B|\nu(\Gamma_0)<\infty$, and hence $N(B \times \Gamma_0)={\rm card}\{i \geq 1; (T_i,X_i,Z_i) \in B \times \Gamma_0\}<\infty$.}
\end{remark}

For any $j \geq 0$, the variable $L_j(B)$ has a compound Poisson distribution with jump intensity measure $|B|\cdot \nu|_{\Gamma_j}$, i.e.
\begin{equation}
\label{ch-funct-Lj(B)}
E[e^{iu L_j(B)}]=\exp\left\{|B|\int_{\Gamma_j}(e^{iuz}-1)\nu(dz) \right\}, \quad u \in \bR.
\end{equation}
It follows that
$E(L_j(B))= |B| \int_{\Gamma_j} z \nu(dz)$ and ${\rm Var}(L_j(B))= |B| \int_{\Gamma_j} z^2 \nu(dz)$  for any $j \geq 0$.
Define
\begin{equation}
\label{def-Y(B)}
Y(B)=\sum_{j \geq 1}[L_j(B)-E(L_j(B))]+L_0(B).
\end{equation}
This sum converges a.s. by Kolmogorov's criterion since $\{L_j(B)-E(L_j(B))\}_{j \geq 1}$ are independent zero-mean random variables with
$\sum_{j \geq 1}{\rm Var}(L_j(B))
<\infty$.

From (\ref{ch-funct-Lj(B)}) and (\ref{def-Y(B)}), it follows that $Y(B)$ is an infinitely divisible random variable with characteristic function:
$$
E(e^{iuY(B)})=\exp \left\{|B|\int_{\bR}(e^{iuz}-1-iu z 1_{\{|z| \leq 1\}})\nu(dz) \right\}, \quad \ u \in \bR.$$
 Hence $E(Y(B))= |B|\int_{\bR}z 1_{\{|z|>1\}}\nu(dz)$ and
${\rm Var}(Y(B))= |B|\int_{\bR}z^2 \nu(dz)$.

In the present article, we assume that
\begin{equation}
\label{v-finite}
v:=\int_{\bR}z^2 \nu(dz)<\infty.
\end{equation}
For any $B \in \cB_{b}(\bR_{+} \times \bR^d)$, we define $$L(B)=Y(B)-E(Y(B))=\sum_{j \geq 0}[L_j(B)-E(L_j(B))].$$
Then $L(B)$ has the characteristic function:
\begin{equation}
\label{ch-func-B}
E(e^{iuL(B)})=\exp \left\{|B|\int_{\bR}(e^{iuz}-1-iu z)\nu(dz) \right\}, \quad \ u \in \bR.
\end{equation}
By Lemma 2.2 of \cite{B13}, the family $\{L(B);B \in \cB_{b}(\bR_{+} \times \bR^d)\}$ is an independently scattered random measure in the sense of \cite{RR89}, with zero mean and covariance: 
\begin{equation}
\label{cov-L-AB}
E[L(A)L(B)]=v|A \cap B|.
\end{equation}

\begin{definition}
{\rm We say that
$L=\{L(B);B \in \cB_{b}(\bR_{+} \times \bR^d)\}$ is a {\em  space-time L\'evy white noise} with jump size intensity $\nu$.}
\end{definition}

Let $\widehat{N}$ be the compensated Poisson measure associated to $N$, i.e. $\widehat{N}(A)=N(A)-E(N(A))$ for any relatively compact set $A$ in $\bR_{+} \times \bR^d \times (\bR \verb2\2 \{0\})$. For any simple function $f=\sum_{i=1}^{n}\alpha_i 1_{A_i}$ on $\bR_{+} \times \bR^d \times (\bR \verb2\2 \{0\})$, we define $$\widehat{N}(f)=\int_{\bR_{+} \times \bR^d \times (\bR \verb2\2 \{0\})} f(t,x,z)\widehat{N}(dt,dx,dz):=\sum_{i=1}^{n}\alpha_i \widehat{N}(A_i).$$
Then $E(\widehat{N}(f))=0$ and $E|\widehat{N}(f)|^2=\int |f(t,x,z)|^2 dtdx \nu(dz)$. By approximation with simple functions, this integral is extended to all functions $f$ with $\int |f(t,x,z)|^2 dtdx \nu(dz)<\infty$.
By (\ref{v-finite}), it follows that for any $B \in \cB_{b}(\bR_{+} \times \bR^d)$,
\begin{equation}
\label{Poisson-rep}
L(B)=\int_{B \times (\bR \verb2\2 \{0\})}z \widehat{N}(dt,dx,dz).
\end{equation}

For any set $B \in \cB_b(\bR^{d+1})$, we define $L(B)=L(B \cap (\bR_{+} \times \bR^d))$ and $L(1_{B})=L(B)$. This definition is extended to simple functions by linearity. For any function $\varphi \in L_{\bC}^2(\bR^{d+1})$, the stochastic integral
$$L(\varphi)=\int_{\bR^{d+1}}\varphi(t,x)L(dt,dx)$$
is defined as a limit in $L^2(\Omega)$, using an approximation by simple functions. Due to (\ref{cov-L-AB}), this integral has the property:
\begin{equation}
\label{cov-L}
E[L(\varphi) \overline{L(\psi)}]=v\int_{\bR^{d+1}}\varphi(t,x) \overline{\psi(t,x)} dtdx.
\end{equation} 

\begin{remark}
\label{spectral-remark}
{\rm The process $\{L(\varphi);\varphi \in \cD_{\bC}(\bR^{d+1})\}$ is a (real) stationary random distribution, in the sense of \cite{ito54}.
Since its covariance $\rho_0:=v\delta_0$ is non-negative definite, there exists a tempered measure $\mu_0$ on $\bR^{d+1}$ (called its spectral measure) such that $\rho_0=\cF \mu_0$ in $\cS_{\bC}'(\bR^{d+1})$. One can easily see that $\mu_0(d\tau,d\xi)= v(2\pi)^{-(d+1)} d\tau d\xi$ ($\tau \in \bR$, $\xi \in \bR^d$), since
$$\rho_0(\varphi)=v \varphi(0,0)=\frac{v}{(2\pi)^{d+1}} \int_{\bR^{d+1}} \cF \varphi(\tau,\xi)d\tau d\xi$$
for any $\varphi \in \cS_{\bC}(\bR^{d+1})$, by the Fourier inversion theorem in $\cS_{\bC}(\bR^{d+1})$. Here $\cF$ denotes the Fourier transform in $(t,x)$.
Therefore, according to Theorem 3 of \cite{yaglom57}, $L(\varphi)$ admits the spectral representation:
$$L(\varphi)=\int_{\bR^{d+1}}\cF \varphi(\tau,\xi)\cM_{0}(d\tau,d\xi) \quad \mbox{for any} \ \varphi \in \cD_{\bC}(\bR^{d+1}),$$
where $\cM_0$ is a symmetric complex random measure with control measure $\mu_0$.  This representation can be extended to all $\varphi \in L_{\bC}^2(\bR^{d+1})$.}
\end{remark}

The next result extends the Poisson representation (\ref{Poisson-rep}) to  $L^2(\bR_{+} \times \bR^{d})$.

\begin{lemma}
\label{properties-L}
For any function $\varphi \in L^2(\bR_{+} \times \bR^{d})$, we have:
\begin{eqnarray*}
(a) & & E(e^{iu L(\varphi)})=\exp \left\{\int_{\bR_{+} \times \bR^d \times \bR} (e^{iuz \varphi(t,x)}-1-iuz \varphi(t,x))dtdx \nu(dz) \right\}, \ u \in \bR \\
(b) & & L(\varphi)=\int_{\bR_{+} \times \bR^d \times (\bR \verb2\2 \{0\})}\varphi(t,x)z \widehat{N}(dt,dx,dz).
\end{eqnarray*}
\end{lemma}

\noindent {\bf Proof:} By Theorem 19.2 of \cite{billingsley95}, there exists a sequence $(\varphi_n)_{n}$ of simple functions on $\bR_{+} \times \bR^d$ such that $|\varphi_n| \leq |\varphi|$ for all $n$ and $(\varphi_n)_n$ converges to $\varphi$ in $L^2(\bR_{+} \times \bR^d)$. By construction, $\{L(\varphi_n)\}_{n}$ converges to $L(\varphi)$ in $L^2(\Omega)$.
By  (\ref{ch-func-B}), for any $n \geq 1$,
$$E(e^{iu L(\varphi_n)})=\exp \left\{\int_{\bR_{+} \times \bR^d \times \bR} (e^{iuz \varphi_n(t,x)}-1-iuz \varphi_n(t,x))dtdx \nu(dz) \right\}$$
Part (a) follows by taking the limit as $n\to \infty$. On the right-hand side, we use the dominated convergence theorem, whose application is justified using the inequality
$|e^{iux}-1-iux| \leq x^2/2$ and (\ref{v-finite}).
By (\ref{Poisson-rep}), for any $n \geq 1$,
$$L(\varphi_n)=\int_{\bR_{+} \times \bR^d \times (\bR \verb2\2 \{0\})}\varphi_n(t,x)z \widehat{N}(dt,dx,dz).$$
Part (b) follows by taking the limit as $n\to \infty$ in $L^2(\Omega)$.
$\Box$

\vspace{3mm}

For any $t>0$ and $\varphi \in L_{\bC}^2(\bR^d)$, we define
$$L_t(\varphi):=L(1_{[0,t]} \varphi).$$
The process $\{L_t(\varphi);t>0, \varphi \in L^2(\bR^d)\}$ is called an ``impulsive cylindrical process on $L^2(\bR^d)$ in Section 7.2 of \cite{PZ07}.
By Lemma \ref{properties-L}.(a), for any $\varphi \in L^2(\bR^d)$,
\begin{equation}
\label{ch-funct-Lt}
E(e^{iu L_t(\varphi)})=\exp\left\{t \int_{\bR^d \times \bR} (e^{iuz \varphi(x)}-1-iuz \varphi(x))dx \nu(dz)  \right\}, \ u \in \bR.
\end{equation}
A similar formula holds for the $L_t(\varphi)-L_s(\varphi)$, and hence
the distribution of $L_t(\varphi)-L_s(\varphi)$ depends only on $t-s$, for any $s<t$.
By Remark \ref{spectral-remark}, this process admits the spectral representation:
\begin{equation}
\label{spectral-repr-L}
L_t(\varphi)=\int_{\bR^d}\cF \varphi(\xi)\cM_t(d\xi) \quad \mbox{for any} \quad \varphi \in L_{\bC}^2(\bR^d),
\end{equation}
where $\cM_{t}(d\xi)= \cF 1_{[0,t]}(\tau)\cM_0(d\tau,d\xi)$ has control measure
$\mu_t(A)=\frac{v}{(2\pi)^{d}}  |A|t$.

\vspace{3mm}

Let
$$\cF_t=\cF_t^N \vee \cN$$ where $\cN=\{F \in \cF; P(F)=0 \ \mbox{or}\ P(F)=1\}$ and $\cF_{t}^N$ is the $\sigma$-field generated by $N([0,s] \times A \times \Gamma)$ for all $s \in [0,t]$, $A \in \cB_b(\bR^d)$ and for all Borel sets $\Gamma \subset \bR \verb2\2 \{0\}$, bounded away from $0$.

\begin{proposition}
(a) For any $t>0$ and $\varphi \in L^2(\bR^d)$, $L_t(\varphi)$ is $\cF_t$-measurable. \\
(b) For any $s<t$ and $\varphi \in L^2(\bR^d)$, $L_t(\varphi)-L_s(\varphi)$
is independent of $\cF_s$.
\end{proposition}

\noindent {\bf Proof:} (a) Without loss of generality, we assume that $\varphi=1_{A}$ with $A \in \cB_{b}(\bR^d)$. (To see this, note that by Theorem 19.2 of \cite{billingsley95}, there exists a sequence $(\varphi_n)_n$ of simple functions such that
$(\varphi_n)_n$ converges to $\varphi$ in $L^2(\bR^d)$. By construction, $\{L_t(\varphi_n)\}_n$ converges to $L_t(\varphi)$ in $L^2(\Omega)$. Hence, $L_{t}(\varphi)=\lim_{k \to \infty}L_{t}(\varphi_{n_k})$ a.s. for a subsequence $(n_k)_k$ and it suffices to prove that $L_t(\varphi_{n_k})$ is $\cF_t$-measurable for any $k$.)
By Lemma \ref{properties-L}.(b),
$$L_t(1_{A})=\int_{[0,t] \times A \times (\bR \verb2\2 \{0\})}z\widehat{N}(ds,dx,dz).$$
Hence, $L_t(1_{A})=\lim_{\varepsilon \to 0}L_{t,\varepsilon}(1_{A})$ in $L^2(\Omega)$, where
$$L_{t,\varepsilon}(1_{A}):=\int_{[0,t] \times A \times \{|z|>\varepsilon\}}z\widehat{N}(ds,dx,dz).$$  It suffices to prove that $L_{t,\varepsilon}(1_{A})$ is  $\cF_t$-measurable, for any $\varepsilon>0$ fixed. 

For this, we approximate the function $h(z)=z1_{\{|z|>\varepsilon\}}$ by a sequence $(h_n)_n$ of simple functions defined as follows: we let $h_n(z)=0$ if $|z| \leq \varepsilon$,
$$h_n(z)=\sum_{k=1}^{n2^n-1}\frac{\varepsilon k}{2^n}1_{(\frac{ \varepsilon k}{2^n},\frac{\varepsilon(k+1)}{2^n}]}(z)+\varepsilon n 1_{(\varepsilon n,\infty)}(z) \quad \mbox{if} \quad z>\varepsilon$$ and $h_n(z)=-h_n(-z)$ is $z<-\varepsilon$.
Then $h_n \to h$ and $|h_n| \leq |h|$ for all $n$. By the dominated convergence theorem, $(h_n)_n$ converges to $h$ in $L^2(\bR,\nu)$.
 It follows that
$$X_n:=\int_{[0,t] \times A \times \{|z|>\varepsilon\}}h_n(z) \widehat{N}(ds,dx,dz)$$
converges in $L^2(\Omega)$ to $L_{t,\varepsilon}(1_{A})$, as $n \to \infty$. Hence, it suffices to prove that $X_n$ is $\cF_t$-measurable for all $n$. This is clear since the sets which appear in the definition of $h_n$ are bounded away from $0$.

(b) Similarly to (a), one can prove that $L_t(\varphi)-L_s(\varphi)$ is $\cF_{s,t}$-measurable, where $\cF_{s,t}$ is the $\sigma$-field generated by $N((a,b] \times A \times \Gamma)$ for all $(a,b] \subset (s,t]$, $A \in \cB_{b}(\bR^d)$ and for all Borel sets $\Gamma \subset \bR \verb2\2 \{0\}$, bounded away from $0$. The conclusion follows since $\cF_{s,t}$ and $\cF_{s}$ are independent. $\Box$

\begin{corollary}
\label{L-Levy}
For any $\varphi \in L^2(\bR^d)$, $\{L_t(\varphi)\}_{t \geq 0}$ is a zero-mean square-integrable L\'evy process with characteristic function (\ref{ch-funct-Lt}). In particular, $\{L_t(\varphi)\}_{t \geq 0}$ is a martingale with respect to $(\cF_t)_{t \geq 0}$.
\end{corollary}

\section{The L\'evy colored noise}
\label{Levy-section}

In this section we introduce an analogue of the spatially homogeneous Gaussian noise considered in \cite{dalang99} for the case of the L\'evy noise. For this, we consider an arbitrary measurable function $h: \bR^d \to \bC$ such that
$|h|^2$ is tempered, i.e.
\begin{equation}
\label{h2-tempered}
\int_{\bR^d}\left(\frac{1}{1+|\xi|^2} \right)^{k}|h(\xi)|^2 d\xi<\infty \quad \mbox{for some} \ k>0.
\end{equation}

\begin{example}
\label{Riesz-ex} (The Riesz kernel)
{\rm Suppose that $h(\xi)=|\xi|^{-\alpha/2}$ for all $\xi \in \bR^d$. Then $|h|^2$ is tempered if and only if $\alpha<d$.
}
\end{example}

Our definition of the colored noise relies on the ``Fourier tranform'' of $L$ in space, which we now define. For any $\phi \in L_{\bC}^2(\bR_{+} \times \bR^d)$, let $\cF^{-1}\phi(t,\cdot)$ be the inverse Fourier transform of $\phi(t,\cdot)$ in $L_{\bC}^2(\bR^d)$, and define
$$\widehat{L}(\phi)=\int_{0}^{\infty}\int_{\bR^d}\phi(t,\xi)\widehat{L}(dt,d\xi):=
\int_{0}^{\infty}\int_{\bR^d}\cF^{-1}\phi(t,\cdot)(x)L(dt,dx).$$
By (\ref{cov-L}) and Plancharel theorem, we see that for any $\phi_1,\phi_2 \in L_{\bC}^2(\bR_{+} \times \bR^d)$,
$$E[\widehat{L}(\phi_1)\overline{\widehat{L}(\phi_2)}]=\frac{v}{(2\pi)^d} \int_{0}^{\infty}\int_{\bR^d}\phi_1(t,\xi)\overline{\phi_2(t,\xi)}d\xi dt.$$

\noindent For any $t>0$ and $\phi \in L_{\bC}^2(\bR^d)$, we let
$\widehat{L}_t(\phi):=\widehat{L}(1_{[0,t]}\phi)=L_t(\cF^{-1}\phi)$. Then
\begin{equation}
\label{cov-hatL}
E[\widehat{L}_t(\phi_1)\overline{\widehat{L}_s(\phi_2)}]=(t \wedge s) \frac{v}{(2\pi)^d}\int_{\bR^d}\phi_1(\xi)\overline{\phi_2(\xi)}d\xi.
\end{equation}
By (\ref{spectral-repr-L}), $\widehat{L}_t(\phi)$ has the spectral representation:
\begin{equation}
\label{spectral-repr-hatL}
\widehat{L}_t(\phi)=\int_{\bR^d}\phi(\xi) \cM_t(d\xi)=:\cM_t(\phi) \quad \mbox{for any} \quad \phi \in L_{\bC}^2(\bR^d).
\end{equation}

\begin{definition}
\label{def-colored-noise}
{\rm For any $t>0$ and $\varphi \in \cS(\bR^d)$, we let:
$$X_t(\varphi):=\int_{0}^{t} \int_{\bR^d} \cF \varphi(\xi)h(\xi)\widehat{L}(ds,d\xi).$$
We say that $\{X_t(\varphi);t \geq 0, \varphi \in \cS(\bR^d)\}$ is a {\em L\'evy colored noise}. }
\end{definition}

\begin{remark}
{\rm $X_t(\varphi)$ is well-defined since the function $\cF \varphi \cdot h$ is in $L_{\bC}^2(\bR^d)$. This follows by (\ref{h2-tempered}), since
$\cF \varphi \in \cS_{\bC}(\bR^d)$ and hence $|\cF \varphi(\xi)|^2 \leq C (1+|\xi|^2)^{-k}$ for all $\xi \in \bR^d$.}
\end{remark}

\begin{remark}
\label{remark-h-Fourier}
{\rm  If there exists a tempered non-negative function $\kappa$ on $\bR^d$ such that
\begin{equation}
\label{h-Fourier-kappa}
h=\cF \kappa \quad \mbox{in} \quad \cS'(\bR^d),
\end{equation}
then $X_t(\varphi)=L_t(\varphi *\kappa)$ for all $\varphi \in \cS(\bR^d)$. (To see this, note that for any $\varphi \in \cS(\bR^d)$, $\varphi *\kappa \in L^2(\bR^d)$ and $\cF(\varphi *\kappa)=\cF \varphi \cdot h$.) In the case of Example \ref{Riesz-ex}, condition (\ref{h-Fourier-kappa}) is satisfied if and only if $\alpha>0$; in this case, $\kappa(x)=c_{\alpha,d}|x|^{-(d-\alpha/2)}$. Condition (\ref{h-Fourier-kappa}) is not needed in the present article.
}
\end{remark}

\begin{remark}
\label{impulsive-remark}
{\rm The L\'evy colored noise is similar to the ``impulsive colored noise'' considered in Section 19.2.2 of \cite{PZ06} in a different framework and more restrictive assumptions (see also  Example 14.26 of \cite{PZ07}). To see this, assume that (\ref{h-Fourier-kappa}) holds. Let $p(dy)=\kappa(y)dy$. By Remark \ref{remark-h-Fourier} and Lemma \ref{properties-L}.(b),
$$X_t(\varphi)=L_t(\varphi * \kappa)=\int_{[0,t] \times \bR^d \times (\bR \verb2\2 \{0\})}  \left(\int_{\bR^d}\varphi(x-y)p(dy) \right)z \widehat{N}(ds,dx,dz).$$
for all $\varphi \in \cS(\bR^d)$. If, in addition, $|h|^2=\cF f$ in $\cS'(\bR^d)$ for some non-negative function $f$ on $\bR^d$, then $1_{A} * \kappa \in L^2(\bR^d)$ for all $A \in \cB_{b}(\bR^d)$ and
$$X_t(A):=L_t(1_{A}*\kappa)=\int_{[0,t] \times \bR^d \times (\bR \verb2\2 \{0\})}  \left(\int_{\bR^d}1_{A}(x-y)p(dy) \right) z \widehat{N}(ds,dx,dz)$$
If condition (19.3) of \cite{PZ06} holds (i.e. $\int_{\Gamma_j}z \nu(dz)=0$ for all $j \geq 0$), then
\begin{eqnarray}
\nonumber
X_t(A)&=& \int_{[0,t] \times \bR^d \times (\bR \verb2\2 \{0\})}  \left(\int_{\bR^d}1_{A}(x-y)p(dy) \right)z N(ds,dx,dz) \\
\label{impulsive-noise}
&=&\sum_{i \geq 1}  1_{[0,t]}(T_i) Z_i \int_{\bR^d} 1_{A}(X_i-y) p(dy)=\sum_{T_i \leq t}Z_i \tilde{p}_{X_i}(A)
\end{eqnarray}
where $\tilde{p}(A)=p(-A)$ and $\tilde{p}_{x}(A)=\tilde{p}(A-x)$ for any $x \in \bR^d$. Relation (\ref{impulsive-noise}) coincides with the representation (19.7) of the impulsive colored noise of \cite{PZ06}.

}
\end{remark}

By (\ref{cov-hatL}), the process $\{X_t(\varphi); t \geq 0, \varphi \in \cS(\bR^d)\}$ has covariance (\ref{cov-X}), where
$$\mu(d\xi)=\frac{v}{(2\pi)^d}|h(\xi)|^2 d\xi.$$

By (\ref{spectral-repr-hatL}), $X_t(\varphi)$ has the spectral representation:
$$X_t(\varphi)=\int_{\bR^d} \cF \varphi (\xi)h(\xi)\cM_t(d\xi), \quad \mbox{for any} \quad \varphi \in \cS(\bR^d).$$

Note that
$$X_t(\varphi)=\widehat{L}_t(\cF \varphi \cdot h)=L_t(\cF^{-1} (\cF \varphi \cdot h)),$$
and hence,  by Corollary \ref{L-Levy}, $\{X_t(\varphi)\}_{t \geq 0}$ is a zero-mean square integrable L\'evy process, and a martingale with respect to $(\cF_t)_{t \geq 0}$. By Theorem 5.4 of \cite{resnick07}, each L\'evy process $\{X_t(\varphi)\}_{t \geq 0}$ has a c\`adl\`ag modification. We work with these modifications.

Recall that, if $M=(M_t)_{t \geq 0}$ is a c\`adl\`ag square-integrable martingale with $M_0=0$, by the Doob-Meyer decomposition, there exists a (unique a.s.) increasing, integrable, right-continuous process $A=(A_t)_{t \geq 0}$ such that $A$ is predictable and $(M_t^2-A_t)_{t \geq 0}$ is a martingale (see e.g. Proposition II.2.1 of \cite{ikeda-watanabe89}). We say that the process $A$ is the {\em predictable variation} of $M$ and we write $A=\langle M \rangle$. (Note that $\langle M \rangle$ may not coincide with the {\em quadratic variation} process $[M]$ defined by $[M]_t =\lim_{n \to \infty} \sum_{j=0}^{k_n-1}(M_{t_{j+1}^n}-M_{t_j^n}) \ \mbox{in} \ L^1(\Omega)$ where $(t_j^n)_{0 \leq j \leq k_n}$ is a partition of $[0,t]$ with $\max_{j}(t_{j+1}^n-t_j^n) \to 0$ as $n \to \infty$. The process $[M]$ is increasing, adapted, c\`adl\`ag, and $(M_t^2-[M]_t)_{t \geq 0}$ is also a martingale.)

For any $\varphi \in \cS(\bR^d)$, the predictable variation of the process $\{X_t(\varphi)\}_{t \geq 0}$ is
$$\langle X_{\cdot}(\varphi)\rangle_t=
t\int_{\bR^d}|\cF \varphi(\xi)|^2 \mu(d\xi).$$

\vspace{3mm}

{\em Stochastic integral with respect to $X$}

\vspace{3mm}

A function $g:\Omega \times \bR_{+} \times \bR^d \to \bR$ is called an {\em elementary process} if
\begin{equation}
\label{elem-g-A}g(\omega,t,x)=Y(\omega)1_{(a,b]}(t)1_{A}(x)
\end{equation}
where $0\leq a <b$, $A \in \cB_b(\bR^d)$ and $Y$ is $\cF_a$-measurable and bounded. We say that $g$ is a {\em smooth elementary process} if it is of the form
\begin{equation}
\label{elem-g}
g(\omega,t,x)=Y(\omega)1_{(a,b]}(t)\psi(x)
\end{equation}
where $0 \leq a<b$, $\psi \in \cD(\bR^d)$ and $Y$ is $\cF_a$-measurable and bounded.
We denote by $\cE$ (respectively $\cE_s$) the set of all linear combinations of elementary processes (respectively smooth elementary processes).

We let $\cP_{\Omega \times \bR_{+}}$ be the predictable $\sigma$-field on $\Omega \times \bR_{+}$ with respect to $(\cF_t)_{t \geq 0}$,
i.e. the $\sigma$-field generated by all linear combinations of processes of the form $X(\omega,t)=Y(\omega) 1_{(a,b]}(t)$, where $0\leq a<b$ and $Y$ is $\cF_a$-measurable and bounded.
A process $\{X(t)\}_{t \geq 0}$ defined on $(\Omega,\cF,P)$ is called {\em predictable} (with respect to $(\cF_t)_{t \geq 0}$) if the map $(\omega,t) \mapsto X(\omega,t)$ is $\cP_{\Omega \times \bR_{+}}$-measurable.

Similarly, we let $\cP_{\Omega \times \bR_{+} \times \bR^d}$ be the {\em predictable} $\sigma$-field on $\Omega \times \bR_{+} \times \bR^d$ with respect to $(\cF_t)_{t \geq 0}$, i.e. the $\sigma$-field generated by all the processes in $\cE$ (or $\cE_s$). Note that
\begin{equation}
\label{pred-inclusion}
\cP_{\Omega \times \bR_{+} \times \bR^d} \subset \cP_{\Omega \times \bR_{+}} \times \cB(\bR^d).
\end{equation}

\begin{definition}
\label{def-pred1}
{\rm A function $g:\Omega \times \bR_{+} \times \bR^d \to \bR$ is called {\em predictable} (with respect to $(\cF_t)_{t \geq 0}$) if it is measurable with respect to $\cP_{\Omega \times \bR_{+} \times \bR^d}$.
}
\end{definition}

For any $g \in \cE_s$ of the form (\ref{elem-g}) and for any $t>0$ , we define
$$(g \cdot X)_t=\int_0^t \int_{\bR^d}g(s,x)X(ds,dx):=Y(X_{t \wedge b}(\psi)-X_{t \wedge a}(\psi)).$$
For any $\varphi \in \cS(\bR^d)$, we let $(g \cdot X)_t(\varphi)=(g \varphi \cdot X)_t$.

The following result follows by classical methods. We omit its proof.

\begin{lemma}
\label{elem-proc-lemma}
For any $g \in \cE_s$ and $\varphi \in \cS(\bR^d)$, $\{(g \cdot X)_t(\varphi)\}_{t \geq 0}$ is a c\`adl\`ag square-integrable martingale with $(g \cdot X)_0(\varphi)=0$,
predictable variation
$$\langle (g \cdot X)_{\cdot}(\varphi) \rangle_{t}=\int_0^t \int_{\bR^d} |\cF(\varphi g(s,\cdot))(\xi))|^2\mu(d\xi)ds$$
and spectral representation
\begin{equation}
\label{spectral-repr-gX}
(g \cdot X)_t(\varphi)=\int_0^t \int_{\bR^d}\cF (\varphi g(s,\cdot))(\xi) h(\xi)\widehat{L}(ds,d\xi).
\end{equation}
In particular, for any $g \in \cE_s$ and $\varphi \in \cS(\bR^d)$,
\begin{equation}
\label{isometry}
E|(g \cdot X)_{t}(\varphi)|^2=E\int_0^t \int_{\bR^d} |\cF(\varphi g(s,\cdot))(\xi))|^2\mu(d\xi)ds.
\end{equation}
\end{lemma}

\begin{remark}
{\rm The term on the right-hand side of (\ref{spectral-repr-gX}) is a {\em stochastic integral with respect to $\widehat{L}$}. This integral is defined as follows. If $g$ is a complex smooth elementary process of the form (\ref{elem-g}) (with $\psi \in \cD_{\bC}(\bR^d)$), we set
$$(g \cdot \widehat{L})_t=Y(\widehat{L}_{t \wedge b}(\psi)-\widehat{L}_{t \wedge a}(\psi)).$$
Then $\{(g \cdot \widehat{L})_t\}_{t \geq 0}$ is a $\bC$-valued square-integrable martingale with variance
$$E|(g \cdot \widehat{L})_t|^2=\frac{v}{(2\pi)^d}E\int_0^t \int_{\bR^d}|g(s,\xi)|^2 d\xi ds.$$
By linearity, this integral is extended to the set $\cE_{s}(\bC)$ of all complex linear combinations of processes of this form. An approximation argument shows that this integral can be extended further to the set of all $\cP_{\Omega \times \bR_{+}} \times \cB(\bR^d)$-measurable functions $g:\Omega \times \bR_{+} \times \bR^d \to \bC$ with
$E\int_0^{\infty} \int_{\bR^d} |g(t,\xi)|^2 d\xi dt<\infty$.}
\end{remark}

Fix $T>0$. As on page 8 of \cite{dalang99}, we introduce the following definition.

\begin{definition}
\label{def-P0}
{\rm Let $\cP_0$ be the completion of $\cE_s$ with respect to $\| \cdot \|_0$, where
$$\|g\|_0^2=E\int_0^T \int_{\bR^d}|\cF g(s,\cdot)(\xi)|^2 \mu(d\xi)ds.$$
}
\end{definition}

The stochastic integral with respect to $X$ can be extended to $\cP_0$ as follows.
The map $g \mapsto \{(g \cdot X)_t\}_{t \in [0,T]}$ is an isometry between $\cE_s$ (endowed with the norm $\|\cdot\|_0$) and the Hilbert space $\cM_2$ of c\`adl\`ag square-integrable $(\cF_t)_t$-martingales $M=(M_t)_{t \in [0,T]}$ with $M_0=0$, equipped with the norm $\|M\|=\{E(M_T^2)\}^{1/2}$. For any $g \in \cP_0$, there exists a sequence $(g_n)_n \subset \cE_s$ such that $\|g_n-g\|_0 \to 0$. By (\ref{isometry}), it follows that $\{(g_n \cdot X)_{t} \}_{t \in [0,T]}, n \geq 1$ is a Cauchy sequence in $\cM_2$. We denote by $g \cdot X=\{(g \cdot X)_t\}_{t \in [0,T]}$ its limit in $\cM_2$ and we write
$$(g \cdot X)_t=\int_0^t \int_{\bR^d}g(s,x)X(ds,dx), \quad t \in [0,T].$$
By construction, for any $g \in \cP_0$, $g \cdot X$ is a c\`adl\`ag square-integrable martingale. In some cases, we can identify its predictable variation (as we will see below).

We proceed now to identify a subset of $\cP_0$, which will be convenient for the study of linear SPDEs with L\'evy colored noise.

\begin{definition}
\label{def-pred2}
{\rm We say that a function $S: \Omega \times [0,T] \to \cS'(\bR^d)$ is {\em predictable} if the map $(\omega,t) \mapsto S(\omega,t,\cdot)(\varphi)$ is $\cP_{\Omega \times \bR_{+}}$-measurable, for any $\varphi \in \cS(\bR^d)$.}
\end{definition}

\begin{remark}
{\rm If $S: \Omega \times [0,T] \to \cS'(\bR^d)$ coincides with a function $g:\Omega \times [0,T] \times \bR^d \to \bR$ (i.e. $S(\omega,t)(\varphi)=\int_{\bR^d} g(\omega,t,x)\varphi(x)dx$ for all $\varphi \in \cS(\bR^d)$) and $g$ is predictable (in the sense of Definition \ref{def-pred1}), then $S$ is predictable (in the sense of Definition \ref{def-pred2}). This follows by (\ref{pred-inclusion}) and Fubini's theorem.
}
\end{remark}

Let $S: \Omega \times [0,T] \to \cS'(\bR^d)$ be a predictable function
such that $\cF S(\omega,t,\cdot)$ is a function for all $(\omega,t)$.
By Lemma 4.2 of \cite{BGP12}, there exists a $\cP_{\Omega \times \bR_+} \times \cB(\bR^d)$-measurable function $\Phi: \Omega \times [0,T] \times \bR^d \to \bC$ such that for all $(\omega,t)$,
\begin{equation}
\label{cF-g}\cF S(\omega,t,\cdot)(\xi)=\Phi(\omega,t,\xi) \quad \mbox{for almost all} \ \xi \in \bR^d.
\end{equation}
Below we will work with $\Phi(\omega,t,\cdot)(\xi)$, but we will write $\cF S(\omega,t,\cdot)(\xi)$.

Let $\overline{\cP}$ be the set of all predictable functions $S: \Omega \times [0,T] \to \cS'(\bR^d)$
such that $\cF S(\omega,t,\cdot)$ is a function for all $(\omega,t)$ and
\begin{equation}
\label{norm-0-g}
\|g\|_0^2:=E\int_{0}^{T}\int_{\bR^d}|\cF S(s,\cdot)(\xi)|^2 \mu(d\xi)dt<\infty.
\end{equation}
This integral is well-defined due to (\ref{cF-g}), the integrand being in fact $\Phi(\omega,t,\cdot)(\xi)$.

The following theorem is the main result of the present article.

\begin{theorem}
\label{main}
Let $S \in \overline{\cP}$ be arbitrary. Then $S \in \cP_0$ and the predictable variation of $S \cdot X$ is
\begin{equation}
\label{pred-var}
\langle S \cdot X \rangle_t=\int_0^t \int_{\bR^d} |\cF S(s,\cdot)(\xi)|^2 \mu(d\xi)ds, \quad t \in [0,T].
\end{equation}
Moreover, $(S \cdot X)_t$ admits the spectral representation:
\begin{equation}
\label{spectral-repr-SX}
(S \cdot X)_t=\int_0^t \int_{\bR^d} \cF S(s,\cdot)(\xi)h(\xi)\widehat{L}(ds,d\xi), \quad t \in [0,T].
\end{equation}
\end{theorem}

\begin{remark}
{\rm The space $\overline{\cP}$ coincides with the space $\Lambda_{X}$ defined on page 20 of \cite{BGP12} (with the measure $F$ replaced by $\mu$, and $\bR_{+}$ replaced by $[0,T]$). The process
$$M=\{M_t(A)=\widehat{L}_t(1_{A}h) ;t \geq 0,A \in \cB_b(\bR^d)\}$$
is an orthogonal martingale measure (as defined in \cite{walsh86}, but with values in $\bC$), and is similar to the process $\{Z_t(A);t \geq 0,A \in \cB_b(\bR^d)\}$ of \cite{BGP12}, except that it is not Gaussian. Relation (\ref{spectral-repr-SX}) can be written in the form:
$$(S \cdot X)_t=\int_0^t \int_{\bR^d} \cF S(s,\cdot)(\xi)M(ds,d\xi), \quad t \in [0,T],$$
and can be viewed as a counterpart of the stochastic integral of \cite{BGP12} (page 21), in the case of the L\'evy noise.}
\end{remark}

\noindent {\bf Proof of Theorem \ref{main}:} We first prove that $S \in \cP_0$. For this, it suffices to prove that for any $\varepsilon>0$, there exists some
 $g_{\varepsilon} \in \cE_s$ such that
\begin{equation}
\label{approx}
E\int_0^T \int_{\bR^d}|\cF g_{\varepsilon}(t,\cdot)(\xi)-\cF S(t,\cdot)(\xi)|^2 \mu(d\xi) dt<\varepsilon.
\end{equation}

Since the function $(\omega,t,\xi) \mapsto \cF S(\omega,s,\cdot)(\xi)$ is $\cP_{\Omega \times \bR_{+}} \times \cB(\bR^d)$-measurable, by applying Theorem 19.2 of \cite{billingsley95} to the real and imaginary part of this function, we infer that there exist some simple $\cP_{\Omega \times \bR_{+}} \times \cB(\bR^d)$-measurable functions $l_n:\Omega \times [0,T] \times \bR^d \to \bC$ such that $l_n(\omega,t,\xi) \to \cF S(\omega,t,\cdot)(\xi)$ and $|l_n(\omega,t,\xi)| \leq |\cF S(\omega,t,\cdot)(\xi)|$ for all $n$. By the dominated convergence theorem, whose application is justified by (\ref{norm-0-g}),
$$E\int_0^T \int_{\bR^d} |l_n(t,\xi)-\cF S(t,\cdot)(\xi)|^2 \mu(d\xi)dt \to 0.$$
This means that for any $\varepsilon>0$ there exists a simple function $l_{\varepsilon}$ such that
\begin{equation}
\label{approx1}
E\int_0^T \int_{\bR^d} |l_{\varepsilon}(t,\xi)-\cF S(t,\cdot)(\xi)|^2 \mu(d\xi)dt <\varepsilon.
\end{equation}
Without loss of generality, we assume that $l_{\varepsilon}(\omega,t,\xi)=1_{F}(\omega,t)1_{A}(\xi)$ for some $F \in \cP_{\Omega \times \bR_{+}},F \subset \Omega \times [0,T]$ and $A \in \cB_{b}(\bR^d)$.

By Lemma 3.8 of \cite{B12}, there exists some $\psi_{\varepsilon} \in \cD_{\bC}(\bR^d)$ such that
$$\int_{\bR^d}|1_{A}(\xi)-\cF \psi_{\varepsilon}(\xi)|^2 \mu(d\xi)<\frac{\varepsilon}{(P \times {\rm Leb})(F)},$$
where ${\rm Leb}$ denotes the Lebesgue measure on $\bR$.

Let $Y(\omega,t)=1_{F}(\omega,t)$. The process $\{Y(t)\}_{t \in [0,T]}$ is predictable (hence, measurable and adapted) and satisfies $E\int_0^T |Y(t)|^2 dt<\infty$
By Lemma II.1.1 of \cite{ikeda-watanabe89}, there exists an elementary process $\{Y_{\varepsilon}(t)\}_{t \in [0,T]}$ on $\Omega \times [0,T]$ such that
$$E\int_0^T |Y_{\varepsilon}(t)-1_F(t)|^2dt<\frac{\varepsilon}{\|\psi_{\varepsilon}\|_0^2},$$
where $\|\psi_{\varepsilon}\|_0^2=\int_{\bR^d} |\cF \psi_{\varepsilon}(\xi)|^2 \mu(d\xi)$. Let $g_{\varepsilon}(\omega,t,x)=Y_{\varepsilon}(\omega,t)\psi_{\varepsilon}(x)$. Then
\begin{eqnarray}
\nonumber
& &  E \int_0^T \int_{\bR^d} |\cF g_{\varepsilon}(t,\cdot)(\xi)-l_{\varepsilon}(t,\xi)|^2 \mu(d\xi)dt=\\
\nonumber
& & E \int_0^T \int_{\bR^d} |Y_{\varepsilon}(t)\cF \psi_{\varepsilon}(\xi)-1_{F}(t)1_{A}(\xi)|^2 \mu(d\xi)dt \leq \\
\nonumber
& & 2 \left(E \int_0^T |Y_{\varepsilon}(t)-1_{F}(t)|^2 dt \int_{\bR^d}|\cF \psi_{\varepsilon}(\xi)|^2 \mu(d\xi)+ \right. \\
\label{approx2}
& & \left. E \int_0^T 1_{F}(t)dt \int_{\bR^d}|\cF \psi_{\varepsilon}(\xi)-1_{A}(\xi)|^2 \mu(d\xi)\right)
< 2(\varepsilon+\varepsilon)=4 \varepsilon.
\end{eqnarray}
Relation (\ref{approx}) follows from (\ref{approx1}) and (\ref{approx2}), since $g_{\varepsilon} \in \cE_s$.

We now prove (\ref{pred-var}). We denote by $A=\{A(t)\}_{t \in [0,T]}$ the process on the right-hand side of (\ref{pred-var}). Clearly, this process is increasing and integrable. Since $A$ is continuous, to prove that it is predictable, it suffices to prove that it is adapted. By definition, $A(t)=\int_0^t V(s)ds$ where
 $$V(s)=\int_{\bR^d}|\cF S(s,\cdot)(\xi)|^2 \mu(d\xi), \quad s \in [0,T].$$  Note that $V=\{V(s)\}_{s \in [0,T]}$ is predictable:
since $(\omega,s,\xi) \mapsto \cF S(\omega,s,\cdot)(\xi)$ is $\cP_{\Omega \times \bR_+} \times \cB(\bR^d)$-measurable, by Fubini's theorem, $(\omega,s) \mapsto V(\omega,s)$ is $\cP_{\Omega \times \bR_+}$-measurable. By Proposition 1.1.12 of \cite{KS91}, $V$ has a progressively measurable modification $\widetilde{V}$. The process $\widetilde{A}(t)=\int_0^t \widetilde{V}(s)ds$ is also progressively measurable, hence adapted. As in the proof of Lemma 3.2.4 of \cite{KS91} (part (c)), one can show that $\widetilde{A}$ is a modification of $A$, i.e. $P(A(t) \not= \widetilde{A}(t))=0$ for all $t$. Since $\cF_t$ contains the $P$-null sets, $A$ is adapted.

Let $M(t)=(S \cdot X)_t$. To prove (\ref{pred-var}), it remains to show that $\{M(t)^2-A(t)\}_{t \in [0,T]}$ is a martingale. This is equivalent to showing that for any $s<t$, and for any $G \in \cF_s$
\begin{equation}
\label{martingale}
E[(M(t)-M(s))^2 1_{G}]=E[(A(t)-A(s)) 1_{G}].
\end{equation}
Since $S \in \cP_0$, there exists a sequence $(g_n)_n \subset \cE_s$ such that \begin{equation}
\label{zero-conv}
E\int_0^T \int_{\bR^d}|\cF g_n(t,\cdot)(\xi)-\cF S(t,\cdot)(\xi)|^2\mu(d\xi)dt \to 0.
\end{equation}
 We denote $M_n(t)=(g_n \cdot X)_t$ and $A_n(t)=\int_0^t \int_{\bR^d}|\cF g_n(s,\cdot)(\xi)|^2 \mu(d\xi)ds$.
By Lemma \ref{elem-proc-lemma}, for any $s<t$, $G \in \cF_s$ and $n \geq 1$,
$$E[(M_n(t)-M_n(s))^2 1_{G}]=E[(A_n(t)-A_n(s))1_{G}].$$
Relation (\ref{martingale}) follows letting $n \to \infty$.
For the left-hand side, we denote $U_n=(M_n(t)-M_n(s))1_{G}$ and $U=(M(t)-M(s))1_G$. By Minkowski inequality,
$$\|U_n-U\|_{L^2(\Omega)} \leq \|(M_n(t)-M(t))1_{G}\|_{L^2(\Omega)}+\|(M_n(s)-M(s))1_{G}\|_{L^2(\Omega)} \to 0,$$
and hence $\|U_n\|_{L^2(\Omega)} \to \|U\|_{L^2(\Omega)}$. For the right-hand side, we use (\ref{zero-conv}).

It remains to prove (\ref{spectral-repr-SX}). Let $(g_n) \subset \cE_s$ be such that (\ref{zero-conv}) holds. By (\ref{spectral-repr-gX}),
$$(g_n \cdot X)_t=\int_0^t \int_{\bR^d} \cF g_n(s,\cdot)(\xi)h(\xi)\widehat{L}(ds,d\xi).$$
Relation  (\ref{spectral-repr-SX}) follows taking the limit as $n \to \infty$ in $L^2(\Omega)$. $\Box$

\vspace{3mm}

{\em Stochastic integral as a martingale measure}

\vspace{3mm}

As in \cite{dalang99}, we suppose now that the following assumption holds:

\vspace{3mm}

\noindent {\em Assumption A.} The Fourier transform of $\mu$ in $\cS'(\bR^d)$ is a non-negative function $f$ on $\bR^d$.

\vspace{3mm}

In this case, for any $\varphi,\psi \in \cS(\bR^d)$,
\begin{equation}
\label{f-mu-rel}
\int_{\bR^d} \int_{\bR^d} \varphi(x) \psi(y) f(x-y)dxdy=\int_{\bR^d} \cF \varphi(\xi) \overline{\cF \psi(\xi)} \mu(d\xi).
\end{equation}
By Lemma 5.6 of \cite{KX09} and polarization, (\ref{f-mu-rel}) holds for any $\varphi,\psi \in L^1(\bR^d)$.

\begin{remark}
{\rm
In the case of Example \ref{Riesz-ex}, Assumption A holds if and only if $\alpha>0$. In this case, $f(x)=c_{\alpha,d}|x|^{-(d-\alpha)}$.}
\end{remark}


For any set $B \in \cB_{b}(\bR^d)$, there exists a sequence $(\varphi_n)_{n} \subset \cD(\bR^d)$ such that $\varphi_n \to 1_{B}$ and ${\rm supp}(\varphi_n) \subset K$ for all $n$ for a compact set $K \subset \bR^d$. The sequence $\{X_t(\varphi_n)\}_n$ is Cauchy in $L^2(\Omega)$ since
\begin{eqnarray*}
\lefteqn{E|X_t(\varphi_n)-X_t(\varphi_m)|^2 =  t\int_{\bR^d}|\cF (\varphi_n-\varphi_m)(\xi)|^2 \mu(d\xi) } \\
& & = t \int_{\bR^d} \int_{\bR^d} (\varphi_n-\varphi_m)(x)(\varphi_n-\varphi_m)(y)f(x-y)dxdy \to 0,
\end{eqnarray*}
by (\ref{cov-X}), (\ref{f-mu-rel}) and the dominated convergence theorem. We denote $X_t(B):=\lim_{n \to \infty}X_t(\varphi_n)$ in $L^2(\Omega)$. Since
$X_t(\varphi_n)=\widehat{L}_t(\cF \varphi_n \cdot h)$ for any $n$,
taking the limit as $n \to \infty$ in $L^2(\Omega)$ we obtain that
$$X_t(B)=\int_0^t \int_{\bR^d}\cF 1_B(\xi)h(\xi)\widehat{L}(ds,d\xi).$$
(The integral on the right-hand side is well-defined since $\cF 1_{B} \cdot h \in L_{\bC}^2(\bR^d)$, due to (\ref{f-mu-rel}).)
Since $$X_t(B)=\widehat{L}_t(\cF 1_{B} \cdot h)=L_t(\cF^{-1}(\cF 1_{B} \cdot h)),$$
by Corollary \ref{L-Levy}, the process $\{X_t(B)\}_{t \geq 0}$ is a zero-mean square-integrable L\'evy process, hence a martingale. This martingale has a c\`adl\`ag modification. We will work with this modification.

It follows that $X=\{X_t(B); t \geq 0, B \in \cB_b(\bR^d)\}$ is a worthy martingale measure (as in \cite{walsh86}), with covariation measure
$$Q_{X}([0,t] \times A \times B):= \langle X_{\cdot}(A),X_{\cdot}(B) \rangle_t=t \int_{A} \int_{B}f(x-y)dxdy$$
and dominating measure $K_X=Q_X$.

If $g$ is an elementary process of the form (\ref{elem-g-A}), and $(\varphi_n)_n \subset \cD(\bR^d)$ is such that $\varphi_n \to 1_{A}$ and ${\rm supp}(\varphi_n) \subset K$ for all $n$ for a compact set $K$, we denote $g_n(\omega,t,x)=Y(\omega)1_{(a,b]}(t)\varphi_n(x)$. By (\ref{cov-X}) and (\ref{f-mu-rel}), $\{(g_n \cdot X)_t\}_{n}$ is a Cauchy sequence in $L^2(\Omega)$. We denote by $(g \cdot X)_t$ its limit.
Since $(g_n \cdot X)_t=Y(X_{t \wedge b}(\varphi_n)-X_{t \wedge a}(\varphi_n))$ for all $n$, taking the limit as $n \to \infty$, we obtain that:
$$(g \cdot X)_t =Y(X_{t \wedge b}(A)-X_{t \wedge a}(A)).$$
For any set $B \in \cB_b(\bR^d)$, we let
$(g \cdot X)_t(B):=(g1_{B} \cdot X)_t$. This definition is extended by linearity to all processes $g \in \cE$.

It follows that for any $g \in \cE$, $g \cdot X=\{(g \cdot X)_t(B); t \geq 0, B \in \cB_{b}(\bR^d)\}$ is also a worthy martingale measure with covariation measure
\begin{eqnarray*}
Q_{g \cdot X}([0,t] \times A \times B)&:=&\langle (g \cdot X)_{\cdot}(A), (g \cdot X)_{\cdot}(B) \rangle_t \\
&=& \int_0^t \int_{A} \int_{B}g(s,x)g(s,y)f(x-y)dxdyds
\end{eqnarray*}
and dominating measure
$$K_{g \cdot X}([0,t] \times A \times B)=\int_0^t \int_{A} \int_{B}|g(s,x)g(s,y)|f(x-y)dxdyds.$$
By approximation, this property continues to hold for all $g \in \cL$, where $\cL$ is the set of predictable functions $g:\Omega \times \bR_{+} \times \bR^d \to \bR$ such that
$$E\int_0^T \int_{B} \int_{B}|g(s,x)g(s,y)|f(x-y)dxdyds<\infty \quad \forall T>0, B \in \cB_b(\bR^d).$$
(Note that $\cE$ is dense in $\cL$.) 


\section{Application to SPDEs}
\label{spde-section}

In this section we consider a linear SPDE driven by the L\'evy colored noise $X$ introduced in Section \ref{Levy-section}.

Let $L$ be a second-order differential operator with constant coefficients. We consider the equation:
\begin{equation}
\label{linear}
Lu(t,x)=\dot{X}(t,x) \quad t>0,x \in \bR^d
\end{equation}
with zero initial conditions.

Let $G$ be the fundamental solution of $Lu=0$. We assume that $G(t,\cdot)$ is a distribution in $\cS'(\bR^d)$ such that its Fourier transform $\cF G(t,\cdot)$ is a function on $\bR^d$ and the map $(t,\xi) \mapsto \cF G(t,\cdot)(\xi)$ is measurable on $\bR_{+} \times \bR^d$.

\begin{example} (Heat equation)
\label{heat-ex}
{\rm Let $Lu=\frac{\partial u}{\partial t}-\frac{1}{2}\Delta u$. Then
$G(t,x)=(2\pi t)^{-d/2} \linebreak \exp\left(-\frac{|x|^2}{2t} \right)$,
where $|\cdot|$ is the Euclidean norm in $\bR^d$. In this case,
$$\cF G(t,\cdot)(\xi)=\exp\left(-\frac{t|\xi|^2}{2} \right)  \quad \mbox{for all} \ \xi \in \bR^d.$$
}
\end{example}

\begin{example} (Wave equation)
\label{wave-ex}
{\rm Let $Lu=\frac{\partial^2 u}{\partial t^2}-\Delta u$. Then  $G(t,\cdot)$ is a function in $L^1(\bR^d)$ if $d=1,2$, a positive measure if $d=3$, and a distribution with rapid decrease if $d \geq 4$. For any $d \geq 1$,
$$\cF G(t,\cdot)(\xi)=\frac{\sin(t|\xi|)}{|\xi|} \quad \mbox{for all} \ \xi \in \bR^d.$$
}
\end{example}

\begin{definition}
{\rm The process $u=\{u(t,x); t \geq 0,x \in \bR^d\}$ defined by
\begin{equation}
\label{def-sol}
u(t,x)=\int_0^t \int_{\bR^d}G(t-s,x-y)X(ds,dy)
\end{equation}
is called a {\em solution} of (\ref{linear}).

}
\end{definition}

By definition, the solution exists if and only if the stochastic integral on the right-hand side of (\ref{def-sol}) is well-defined, i.e. $G(t-\cdot,x-\cdot)1_{[0,t]} \in \cP_0$, where $\cP_0$ is the space given by Definition \ref{def-P0} with $T$ replaced by $t$.

The following result is proved similarly to Theorem \ref{main}. We omit the details.
Note that Assumption A is not required for this result.

\begin{theorem}
\label{linear-th}
Equation (\ref{linear}) has a solution $u$ if and only if
\begin{equation}
\label{I-finite}
I_t:=\int_{0}^{t} \int_{\bR^d}|\cF G(s,\cdot)(\xi)|^2 \mu(d\xi)ds<\infty, \quad \forall t>0.
\end{equation}
In this case, $E|u(t,x)|^2=I_t$ and $u(t,x)$ admits the spectral representation:
$$u(t,x)=\int_0^t \int_{\bR^d} e^{-i \xi \cdot x} \overline{\cF G(t-s,\cdot)(\xi)}h(\xi)\widehat{L}(ds,d\xi).$$
\end{theorem}

\begin{remark}
{\rm In Examples \ref{heat-ex} and \ref{wave-ex}, condition (\ref{I-finite}) holds if and only if
\begin{equation}
\label{cond-mu}
\int_{\bR^d}\frac{1}{1+|\xi|^2}\mu(d\xi)<\infty.
\end{equation}
When $h(\xi)=|\xi|^{-\alpha/2}$ for some $\alpha<d$ (Example \ref{Riesz-ex}), (\ref{cond-mu}) holds if and only if $\alpha>d-2$. When $d \geq 2$, this implies that $\alpha \in (0,d)$ (hence Assumption A holds). But when $d=1$, Theorem \ref{linear-th} is valid for any $\alpha \in (-1,1)$. }
\end{remark}

\noindent \footnotesize{{\em Acknowledgement.} The author is grateful to Robert Dalang for suggesting this problem.

\normalsize{

\end{document}